\documentclass[12pt]{amsart}
\usepackage{amssymb}
\usepackage{amsthm}
\usepackage{euscript}

\theoremstyle{plain}

\newtheorem{thm}[subsection]{Theorem}
\newtheorem{prop}[subsection]{Proposition}

\newtheorem{cor}[subsection]{Corollary}

\begin{document}
\title[Extensions of the fundamental theorem of algebra ]{Extensions of the fundamental theorem of algebra }
\author{Bamdad R. Yahaghi}

\address{Department of Mathematics, Faculty of Sciences, Golestan University, Gorgan 19395-5746, Iran}
\email{ bamdad5@hotmail.com, bbaammddaadd55@gmail.com}



\bibliographystyle{plain}

\begin{abstract}
In this paper motivated by the celebrated fundamental theorem of algebra and its standard proof utilizing Liouville's Theorem, we prove the fundamental theorem of algebra type results for both commutative and noncommutative polynomials in the setting of left (resp. right) alternative topological  complex algebras whose topological duals separates their elements and that of such real algebras whose centers contain certain copies of complex numbers. An application of one of the main results of the paper is the existence of eigenvalues for matrices with entries from arbitrary finite-dimensional complex algebras. We also prove the existence of right eigenvalues for  matrices with entries from finite-dimensional associative real algebras that contain  copies of the complex numbers.  

\end{abstract}

\maketitle 

\bigskip

\begin{section}
{\bf Introduction}
\end{section}

\bigskip

The celebrated fundamental theorem of algebra (FTA) needs no introduction.  It is safe to say that we all were exposed to the FTA, perhaps the real version of it, in high school. From then on, with the fundamental theorem at our disposal, we did various interesting problems and theorems such as {\it any polynomial with real or complex coefficients that is  nonnegative on the real line is the sum of squares of two real polynomials} and  Lucas' theorem which asserts {\it the roots of the derivative of a given nonconstant complex polynomial lie in the convex hull of those of the given polynomial.} All that being said, we refer the reader to \cite[Chapter 4]{EHHKMNPR} and \cite{FR}  for a detailed introduction as well as a comprehensive account of the fundamental theorem and its history.  For a fundamental theorem of algebra for polynomial
equations over real composition algebras, see \cite{W}.

In this note, we consider arbitrary   real or complex unital algebras of the following kinds: finite-dimensional algebras and left (resp. right) alternative topological algebras whose topological duals separate their elements. First, we prove a Fundamental Theorem of Algebra type result  for both commutative and noncommutative polynomials with coefficients from  left (resp. right) alternative topological  algebras whose topological duals separate their elements.   We use this result to prove the existence of eigenvalues for matrices with entries from arbitrary finite-dimensional complex unital algebras as well as real unital algebras that contain  certain copies of the complex numbers.  Finally, we prove the existence of right eigenvalues for matrices with entries in  finite-dimensional associative real algebras that contain  copies of the complex numbers.

Let's  get down to work by setting the stage. A vector space  $\mathbb{A}$ over reals (resp. complex numbers) together with a multiplication coming from a bilinear form on  $\mathbb{A}$  is said to be a {\it real (resp. complex) algebra}. The algebra  $\mathbb{A}$ is said to be {\it unital} or {\it to have an identity element}  if its multiplication has an identity element, denoted by $1$. The identity element of the addition operation of the algebra is denoted by $0$. The algebra  $\mathbb{A}$ is called {\it associative (resp. commutative)}  if its multiplication is associative (resp. commutative). Throughout, by an algebra we mean an arbitrary real or complex nonzero algebra not necessarily associative or commutative.  Also, as is usual, the symbol $\mathbb{F}$ stands for $\mathbb{R}$ or $\mathbb{C}$.  A nonzero element $ a \in  \mathbb{A}$ is said to be {\it invertible} if there exists a unique element of the algebra, denoted by $a^{-1}$, satisfying the relations $ a a^{-1} = a^{-1} a = 1$. The symbol  $\mathbb{A}^{-1}$ is used to denote the set of all invertible elements of the algebra  $\mathbb{A}$. 

For $ a, b , c \in \mathbb{A}$, by definition 
$$ [a, b ] := ab - ba, \ \ \ [a, b , c] :=(ab)c - a(b c),$$
are called {\it the commutator of the elements} $a, b$  and  {\it the associator of the  elements} $ a, b , c$, respectively. It is plain that the commutator and the associator are bilinear and trilinear functions, respectively.  An $ a \in \mathbb{A}$ is said to be a {\it central element} if $ [ a, \mathbb{A} ] :=  \{ [ a , b ] : b \in \mathbb{A} \} = \{0\}$. An $ a \in \mathbb{A}$ is said to be a {\it nuclear element} if $ [ a, \mathbb{A}, \mathbb{A} ] = [  \mathbb{A}, a,  \mathbb{A} ]  =  [  \mathbb{A},  \mathbb{A}, a ]  =  \{0\}$. The symbols ${\rm Z}(\mathbb{A}) $ (resp. ${\rm N}(\mathbb{A} ) $) are used to denote the {\it centre (resp. nucleus)} of $  \mathbb{A}$, namely, the set of all central (resp. nuclear) elements of the algebra   $  \mathbb{A}$. Note that the subalgebra generated by any nuclear element of an arbitrary algebra is both commutative and associative.

An algebra  $\mathbb{A}$ is said to be {\it left (resp. right) alternative} if  $ a(ab) = (aa)b $ (resp. $ a(bb) = (ab)b$) for all $ a, b \in  \mathbb{A}$, which are respectively called the left and the right alternative identities. The algebra  $\mathbb{A}$ is said to be {\it alternative} if it is both left and right alternative;   $\mathbb{A}$ is said to be {\it one-sided alternative} if it is  left or right alternative. By a theorem of E. Artin, \cite[Theorem 3.1]{Sch}, an algebra  $\mathbb{A}$ is alternative if and only if for every $a, b \in  \mathbb{A}$, the subalgebra generated by the elements $ a$ and $b$ is associative.   An algebra  $\mathbb{A}$ is said to be {\it power-associative} if, for every $a \in  \mathbb{A}$, the subalgebra generated by the element $ a$ is associative. By a theorem of Albert, \cite[Lemma 1]{A2}, one-sided alternative algebras over a general field $F$  are power-associative whenever  the characteristic of $F$ is not $2$.

The following identities known as {\it the Moufang's identities} (see  \cite[Page 28]{Sch} or \cite[Lemma 2.3.60]{GP}) hold on  any alternative algebra $\mathbb A$
$$ (ab a)c = a \big( b (a c)  \big)     , \ \ \   (ab) (ca) = a (bc) a  , \  \ \       a(bc b) = \big( (ab) c \big) b.   $$
Also the following identities,  which we call {\it the Skornyakov's Identities} (see \cite[Pages 178-179]{Sk}), hold on  any left (resp. right) alternative algebra $\mathbb A$ over a general field $F$ whose characteristic  is not $2$
$$ (a (b a))c = a \big( b (a c) \big)     , \ \ \  [a, b, c] (cb) = \big([a, b, c] b\big) c,   $$
$$ \Big ({\rm resp.} \       a((bc) b) = \big( (ab) c \big) b , \ \ \   (ba) [a, b, c] = a \big( b  [a, b, c] \big). \Big) $$

 Note that if an algebra  $\mathbb{A}$ is associative, then the uniqueness of the inverse element is a redundant hypothesis in the definition of the invertible elements of  $\mathbb{A}$. It turns out that the uniqueness in the definition of invertible elements of an algebra  $\mathbb{A}$ is a redundant hypothesis whenever the algebra  $\mathbb{A}$ is left (resp. right) alternative, see \cite[Proposition 2.3(iv)]{Y}.  

By a {\it topological real (resp. complex) algebra}, we mean a   real (resp. complex) algebra  $\mathbb{A}$ together with a Hausdorff topology with the property that the operations addition and scalar product of the algebra $\mathbb{A}$  are continuous, the multiplication is separately continuous, and that the inversion function on $\mathbb{A}^{-1}$  is continuous.  For a topological algebra $\mathbb{A}$ over $\mathbb{F}$,  we use the symbols $\mathcal{L}_c(\mathbb{A})$  (resp. $ \mathbb{A}^* := \mathcal{L}_c(\mathbb{A}, \mathbb{F}) $) to denote the set of all continuous $\mathbb{F}$-linear operators (resp. functionals)  from  $\mathbb{A}$  into $\mathbb{A}$ (resp.  $\mathbb{F}$);  $ \mathbb{A}^* $ is called the {\it (topological) dual} of the algebra $\mathbb{A}$. The topological algebra  $\mathbb{A}$ is said to be {\it locally convex} if its topology is generated by a family   $\{s_j\}_{j \in J}$ of real (resp. complex) vector space seminorms.  By an {\it lmc (locally multiplicatively convex) algebra}, we mean an algebra equipped with a locally convex vector space topology whose seminorms are algebra seminorms. One can show that lmc algebras are topological algebras.  It is worth mentioning that by the Hahn-Banach Theorem, \cite[Theorem 3.4(b)]{Ru}, the topological duals of locally convex vector spaces separate the elements of the spaces.  An algebra norm $||.||$  of a unital algebra  $\mathbb{A}$ is said to be unital if $||1||= 1$,  where the first $1$ denotes the identity  element of the algebra  $\mathbb{A}$.

Let $R$ be a commutative ring with identity.  As is usual, the symbol $R[x]$ stands for the ring of all polynomials in the indeterminate $x$ with coefficients in the ring $R$. Let $ f= f_0 + f_1 x + \cdots + f_n x^n \in R[x]$ be of degree $n \in \mathbb{N}_0 := \mathbb{N} \cup \{0\}$, i.e., $ f_n \not= 0$, where $R$  is  a commutative ring with identity; the coefficient $f_n \in R\setminus \{0\}$ is called the leading coefficient of the polynomial $f$. We say that an element $ r \in R$ is a singular element for the polynomial $f  $ if $ f(r):= f_0 + f_1 r + \cdots + f_n r^n \notin R^{-1}$.  Let $R$ be a ring with identity, not necessarily commutative nor associative. An expression of the form $ f_0 x f_1 x f_2 \cdots x f_n$ with  $n \in \mathbb{N}_0$ and $ f_i \in R \setminus \{ 0\}$ ($0 \leq i \leq n$) is said to be a noncommutative monomial of degree $n$ with coefficients $ f_i$   in the indeterminate $x$. Note that if the ring is nonassociative, proper parentheses must be inserted in the expression $ f_0 x f_1 x f_2 \cdots x f_n$ to make it sensible.  By a noncommutative polynomial in the indeterminate $x$ and with coefficients in the ring $R$, we mean a finite sum of noncommutative monomials with coefficients in $R$, each of which is called a monomial summand of the noncommutative polynomial. The sum of all noncommutative monomials of the greatest degree in a noncommutative polynomial is called the leading polynomial part of the noncommutative polynomial. For instance, $f_0 +  f_1 (x f_1')  +   (f_1 x) f_1'+ f_1'' x + xf'''_1 $, where the coefficients come from $R$,  is an example of a noncommutative polynomial of degree at most $1$ with coefficients in $R$ whose leading noncommutative polynomial part is  $  f_1 (x f_1')  +   (f_1 x) f_1'+ f_1'' x + xf'''_1 $. The notion of singular elements of noncommutative polynomials with coefficients in a unital ring $R$ can be defined similarly. 

Let $R$ be a  ring and $ A \in M_n (R)$, the set of all $n \times n$ matrices with entries from $R$. An element  $ \lambda \in R$ is said to be a left (resp. right) eigenvalue of the matrix $A$ if there is a nonzero   $ n \times 1$ column matrix $ X \in R^n := M_{n \times 1} (R)$ such that $ AX = \lambda X $ (resp.  $ AX =X  \lambda $). An element  $ \lambda \in R$ is said to be an eigenvalue of the matrix $A$ if there is a nonzero   $ n \times 1$ column matrix $ X \in R^n := M_{n \times 1} (R)$ such that $ AX = \lambda X =X  \lambda $.

\bigskip

The first assertion of part (i) of the following theorem is a slight extension of \cite[Proposition 1.1.7]{GP} and \cite[Theorem 4]{A}.  Part (ii) of the theorem might be of independent interest. For a thorough treatment of the theory of (complete) normed algebras, we refer the reader to the classical references   \cite{Ri} and \cite{BD} and to the more recent excellent reference \cite{GP}. 

\bigskip

\begin{thm} \label{1.1} 
 {\rm (i)}  Let $\mathbb{A}$ be a finite-dimensional real or complex algebra and $||.||$ be a vector space norm on $\mathbb{A}$. Then, there exists an $M > 0$ such that $M||.||$ is an algebra norm on $\mathbb{A}$.  Therefore, every finite-dimensional real or complex algebra can be normed. Moreover, if a real or complex algebra is finite-dimensional and unital, it can be equipped  with a unital norm.

 {\rm (ii)}   Let $\mathbb{A}$ be a finite-dimensional real or complex algebra and $s$ be a vector space semi-norm on $\mathbb{A}$. Then, there exists an $M > 0$ such that $Ms$ is an algebra semi-norm on $\mathbb{A}$ if and only if the kernel of  $s$ is a two-sided ideal of $\mathbb{A}$. 
\end{thm}  

\bigskip

\noindent {\bf Proof.} (i)  The first assertion implies the second assertion because every real or complex vector space can be normed. Now, choose an ordered  basis $\mathcal{B}= (u_i)_{i=1}^n$ for  $\mathbb{A}$. Clearly,  $|| .||_1= ||.||_{\mathcal{B}1} :  \mathbb{A}   \longrightarrow \mathbb{R}$ defined by $ ||a||_1 =  ||a||_{\mathcal{B}1} =\sum_{i=1}^n |a_i|$, where $ a = \sum_{i=1}^n a_i u_i$ and $ a_i \in \mathbb{F}$ 
($1 \leq i \leq n$),  is a vector space norm on  $\mathbb{A}$.  Since $\mathbb{A}$ is finite-dimensional, the two norms $||.||$ and $|| .||_1$ on  $\mathbb{A}$ are equivalent. Thus, there are $ c_1 , c_2 > 0$ such that $ c_1 ||.|| \leq ||.||_1 \leq c_2 ||.||$. 
Set $M_0:= \max_{1 \leq i , j  \leq n} ||u_i  u_j||$. 
Given $ a = \sum_{i=1}^n a_i u_i$ and $ b = \sum_{j=1}^n b_j u_j$ in  $\mathbb{A}$ with $ a_i , b_i \in  \mathbb{F} $
($1 \leq i \leq n$),  we can write 
\begin{eqnarray*}
  ||ab||  & = &  \Big\|  \sum_{i, j=1}^n a_i b_j u_i u_j   \Big\|  \\
  & \leq &     \sum_{i, j=1}^n |a_i|  | b_j| ||u_i u_j||, \\
  & \leq &  M_0 \sum_{i, j=1}^n |a_i |  | b_j| =  M_0 ||a||_1 ||b||_1 , \\
   & \leq & c_2^2 M_0 ||a|| ||b|| 
\end{eqnarray*}
implying that 
$|| ab || \leq c_2^2 M_0 ||a|| ||b||$, equivalently,  
$$ M|| ab || \leq \Big(M ||a||\Big) \Big( M||b|| \Big), $$ 
 where $M:= c_2^2 M_0$. Therefore, $ ||.||' :  \mathbb{A}   \longrightarrow \mathbb{R}$ defined by $ ||a||' = M ||a||$ is an algebra norm on  $\mathbb{A}$, as desired. 

Next, let the algebra  $\mathbb{A}$ be unital. Use the first part to equip  $\mathbb{A}$ with an  algebra norm $||.||$ on  $\mathbb{A}$. Define $||.||' :  \mathbb{A}   \longrightarrow \mathbb{R}$
by  $||a||' := \inf_{\lambda \in \mathbb{F}} ( |\lambda| +  ||a - \lambda||)$. As shown in   \cite[Proposition 1.1.111]{GP}, $||.||' $ is a unital algebra norm on  $\mathbb{A}$, which is what we want. 

(ii) The necessary part of the assertion is obvious and it holds even when the ground algebra is infinite-dimensional. As for the sufficiency part,  it holds even when the algebra $\mathbb{A}$ is infinite-dimensional. But the two-sided ideal $ \mathcal{K}$, the kernel of  $s$, must be of a finite codimension, i.e.,  $\mathbb{B} := \frac{\mathbb{A}}{\mathcal{K}}$ must be finite-dimensional. With this hypothesis at our disposal, define $s_\mathbb{B} : \mathbb{B} \longrightarrow \mathbb{R}$ by  $s_\mathbb{B} ( a + \mathcal{K}) := s(a)$. Since $ \mathcal{K} := \ker s$, we get that $s_\mathbb{B}$ is well-defined and indeed a vector space norm on the finite-dimensional algebra $\mathbb{B}$. It thus follows from (i) that there is an $M > 0$ such that $Ms_\mathbb{B} $ is an algebra norm on $\mathbb{B}$. Now, given arbitrary $ x, y \in  \mathbb{A}$, we can write
\begin{eqnarray*}
Ms(xy) & = &  M s_\mathbb{B} ( xy + \mathcal{K}), \\
& =& M s_\mathbb{B} \big( (x + \mathcal{K} ) ( y +  \mathcal{K} ) \big), \\
& \leq & M s_\mathbb{B} \big( x + \mathcal{K}  \big)  M s_\mathbb{B} \big( y +  \mathcal{K}  \big), \\
& = & Ms(x)Ms(y).
\end{eqnarray*}
In other words,  the vector space semi-norm $Ms $ is an algebra semi-norm on $\mathbb{A}$, as desired.  
\hfill \qed

\bigskip

\noindent {\bf Remark.} Part (i) of the theorem can be deduced from the following proposition, which is known by the experts.  {\it Let $\mathbb{A}$ be a real or complex algebra and $||.||$ be a vector space norm on $\mathbb{A}$ with respect to which the multiplication of $\mathbb{A}$ is continuous at $(0, 0)$. Then, there exists an $M > 0$ such that $M||.||$ is an algebra norm on $\mathbb{A}$. } To prove this, just note that the continuity of the multiplication at $(0, 0)$ implies that there exists an $r > 0$ such that $\| a b \| < 1$ whenever $ \|a\| <r $ and $ \|b\| <r$. This, in turn, implies that $ \| ab\| \leq r^2 \| a\| \|b\|$ for all $ a, b \in  \mathbb{A}$, from which we get that $ r^2 \|.\|$ is  an algebra norm on $\mathbb{A}$, as desired.  Indeed, the aforementioned proposition, in turn, can be deduced from the following fact, known by the experts, about bilinear functions on NLS'. {\it Let $X_1$,  $X_2$, and $Y$ be real (resp. complex) NLS' and $T :  X_1 \times X_2 \longrightarrow Y$ be a bilinear function, where $X_1 \times X_2$ is endowed with the product topology. Then, $T $ is continuous if and only if it is continuous  at $(0, 0)$, if and only if $\|T\|  < +\infty $, where  
$$\|T\| := \sup \Big\{ \|T(x_1, x_2) \|_Y :  \|x_i\|_{X_i} \leq 1, \ i = 1, 2 \Big\} =$$
$$ \inf \Big\{ c > 0:   \|T(x_1, x_2) \|_Y \leq c   \|x_1\|_{X_1} \|x_2\|_{X_2}  \ \forall x_i \in  X_i,  i = 1, 2 \Big\}. $$}

\bigskip

\begin{cor} \label{1.2} 
Let $\mathbb{A}$ be a finite-dimensional real or complex algebra and $s$ be a vector space semi-norm on $\mathbb{A}$. Then, there exist an $M > 0$ and algebra norm $||.||$ on $\mathbb{A}$ such that $Ms$ is dominated by $||.||$, i.e.,  $Ms \leq ||.||$.
\end{cor}  

\bigskip

\noindent {\bf Proof.} Let $ \mathcal{K} := \ker s := \{ a \in \mathbb{A}: s(a) = 0\}$, which is a subspace of $\mathbb{A}$. Let $\mathcal{M}$ be a complementary subspace of $\mathbb{A}$, i.e., $ \mathcal{K} + \mathcal{M} = \mathbb{A}$ and $ \mathcal{K} \cap \mathcal{M} =\{ 0\}$. Let $P_ \mathcal{K}$ be the projection onto $\mathcal{K}$ along $ \mathcal{M} $. Choose a vector space norm $n$ on the subspace $\mathcal{K}$. It is plain that $ \| .\|_1 : \mathbb{A} \longrightarrow \mathbb{R}$ defined by $ \| .\|_1 := s + nP_\mathcal{K}$  is a vector space norm on  $\mathbb{A}$. The assertion now follows from the preceding theorem. 
\hfill \qed

\bigskip

\begin{section}
{\bf Main results}
\end{section}

\bigskip

We start off with a useful proposition.

\bigskip 

\begin{prop} \label{2.1} 
   {\rm (i)} Let $( \mathbb{A}, ||.||)$ be a  normed alternative  unital real algebra. If the algebra  $\mathbb{A}$  contains a copy of complex numbers, say,  $ \mathbb{C}_I := \{ a + b I: a, b \in \mathbb{R}\}$ with $I \in  \mathbb{A}$ and $ I^2 = -1$,  then there exists an algebra norm $ ||.||' :  \mathbb{A}  \longrightarrow \mathbb{R}$ on  $\mathbb{A}$ equivalent to  $ ||.||$ such that  $|| z a w||' = 
 |z| ||a||' |w|$ for all $ z, w \in \mathbb{C}_I $ and $ a \in  \mathbb{A}$. 

   {\rm (ii)} Let $( \mathbb{A}, ||.||)$ be a  normed left (resp. right) alternative  unital real algebra. If the algebra  $\mathbb{A}$  contains a copy of complex numbers, say,  $ \mathbb{C}_I := \{ a + b I: a, b \in \mathbb{R}\}$ with $I \in  \mathbb{A}$ and $ I^2 = -1$,  then there exists an algebra norm $ ||.||' :  \mathbb{A}  \longrightarrow \mathbb{R}$ on  $\mathbb{A}$ equivalent to  $ ||.||$ such that 
    $|| z a ||' = |z| ||a||'$  (resp.  $|| a w||' =||a||' |w|$)  for all $ z, w \in \mathbb{C}_I $ and $ a \in  \mathbb{A}$. 

\end{prop}  

\bigskip

\noindent {\bf Proof.} (i)  Let $c:= \sup_{t \in \mathbb{R}} \|e^{I t}\| $.      We show that $\| .\|' :  \mathbb{A} \longrightarrow \mathbb{R}$ defined by 
$$\|a\|' := c \sup \{ \|e^{I t_1} a e^{ I t_2}\| :  t_i \in \mathbb{R},  i = 1, 2 \} $$
is  indeed an algebra norm satisfying      $ c \| .\| \leq      \| .\|' \leq c^3 \| .\|$ and 
 $|| z a w||' = 
 |z| ||a||' |w|$ for all $ z, w \in \mathbb{C}_I $ and $ a \in  \mathbb{A}$. That $|| .||'$ is a vector space norm and satisfies the desired properties is clear. It remains to be shown that 
  $|| .||'$ is an algebra norm.  To this end, because of the alternativity of the algebra and the middle Moufang's identity, for given $ a, b \in  \mathbb{A}$, we can write 
 \begin{eqnarray*}
 ||a b||' & =&  c \sup \left\{ ||e^{I t_1} (ab) e^{ I t_2}|| :  t_i \in \mathbb{R},  i = 1, 2 \right\},\\
 & =& c \sup \left\{ \Big\|\big(e^{I t_1} (ab) e^{ I t_1}\big) e^{I(t_2 - t_1)} \Big\| :  t_i \in \mathbb{R},  i = 1, 2 \right\},\\
 &\leq &  c^2 \sup \left\{ \Big\|\big(e^{I t} a \big) \big( b e^{ I t}\big) \Big\| :  t \in \mathbb{R} \right\}, \\
  &\leq &  c^2  \sup_{t \in \mathbb{R}}  \big\| e^{I t} a \big\|      \sup_{t \in \mathbb{R}} \big\| b e^{ I t} \big\| , \\
  &\leq & \| a\|' \|b\|'.
 \end{eqnarray*}

 (ii) It suffices to prove the assertion for left alternative algebras.   Let $c:= \sup_{t \in \mathbb{R}} \|e^{I t}\| $.      We show that $\| .\|' :  \mathbb{A} \longrightarrow \mathbb{R}$ defined by 
$$\|a\|' := c \sup \{ \| \big(e^{I t_1} a\big)  e^{ I t_2}\| :  t_i \in \mathbb{R},  i = 1, 2 \} $$
is  indeed an algebra norm satisfying      $ c \| .\| \leq      \| .\|' \leq c^3 \| .\|$ and 
 $|| z a ||' = 
 |z| ||a||' $ for all $ z \in \mathbb{C}_I $ and $ a \in  \mathbb{A}$. That $|| .||'$ is a vector space norm and satisfies   $ c \| .\| \leq      \| .\|' \leq c^3 \| .\|$  is clear. It remains to be shown that  $|| z a ||' = 
 |z| ||a||' $ for all $ z \in \mathbb{C}_I $ and $ a \in  \mathbb{A}$ and that $\| .\|' $ is submultiplicative, which means it is indeed an algebra norm on $\mathbb{A}$. To this end, first, fix  $ z \in \mathbb{C}_I $ and $ a \in  \mathbb{A}$. It follows that there exists a $ \theta \in \mathbb{R}$ such that $ z = |z| e^{I \theta}$. So we can write 
 \begin{eqnarray*}
 \|za\|'  & =&  c \sup \{ \| \big(e^{I t_1} (za)\big)  e^{ I t_2}\| :  t_i \in \mathbb{R},  i = 1, 2 \} \\
  & = &  |z|  c \sup \{ \| \big(e^{I t_1} (e^{I \theta}a)\big)  e^{ I t_2}\| :  t_i \in \mathbb{R},  i = 1, 2 \} \\
   & = & |z|  c \sup \{ \| \big(e^{I (t_1 + \theta)} a \big)  e^{ I t_2}\| :  t_i \in \mathbb{R},  i = 1, 2 \} \\
    & = & |z| ||a||',
    \end{eqnarray*}
 as desired. Next, fixing arbitrary $ a, b \in  \mathbb{A}$, utilizing Skorniyakov's left identity, we have 
  \begin{eqnarray*}
 \|ab\|'  & =& c \sup \{ \| \big(e^{I t_1} (ab)\big)  e^{ I t_2}\| :  t_i \in \mathbb{R},  i = 1, 2 \} \\
  & =&   c \sup \{ \| e^{I t_1}\big(a (b  e^{ I t_2}) \big)\| :  t_i \in \mathbb{R},  i = 1, 2 \}\\
    & \leq & c \sup \{ \| e^{I t_1}\| \|a\| \| b  e^{ I t_2}\| : t_i \in \mathbb{R},  i = 1, 2 \}\\
     & \leq & \| a\|' \|b\|',
     \end{eqnarray*}
 which is what we wanted. This completes the proof.   
\hfill \qed 

\bigskip

The following theorem can be thought of as an extension of the Fundamental Theorem of Algebra to  commutative and  left or  right alternative topological real (resp. complex) algebras whose multiplications are continuous and whose topological duals separate their elements. In the proofs of the next two theorems, we have made use of the identities 
 \begin{eqnarray*}
  a^{-1} - b^{-1} & = &   a^{-1}\big[  ( b-a) b^{-1} \big], \\
 \Big( {\rm resp.} \   a^{-1} - b^{-1} & = &  \big[ a^{-1}  ( b-a)  \big] b^{-1}, \Big)
   \end{eqnarray*}
which, as pointed out preceding \cite[Theorem 2.4]{Y}, holds  on any left (resp. right) alternative algebra. The aforementioned identity as explained in \cite{Y} prepares the ground to develop the spectral theory in the setting of normed, more generally topological, alternative complex (and hence real) algebras.

\bigskip 

\begin{thm} \label{2.2} 
  {\rm (i)} 
Let $\mathbb{A}$  be a  commutative and left (resp. right) alternative topological unital real algebra whose multiplication is continuous and whose topological dual separates its points, e.g., a locally convex real algebra whose multiplication is jointly continuous. Asuume further that the algebra $\mathbb{A}$ contains a copy of the complex numbers, say, $ \mathbb{C}_I := \{ a + b I: a, b \in \mathbb{R}\}$, where $I \in \mathbb{A}$ is such that $ I^2 = -1$. Then, every nonconstant polynomial with coefficients in $\mathbb{A}$ and with an invertible leading coefficient has a singular element in  $ \mathbb{C}_I $.  

  {\rm (ii)}  Let $\mathbb{A}$  be a  commutative and  left (resp. right) alternative topological unital complex  algebra whose multiplication is continuous and whose topological dual separates its points, e.g., a locally convex complex algebra whose multiplication is jointly continuous.  Then, every nonconstant polynomial with coefficients in $\mathbb{A}$ and with an invertible leading coefficient has a singular element in  $ \mathbb{C} 1 $.  

\end{thm}  

\bigskip

\noindent {\bf Proof.} (i) We prove the assertion for left alternative algebras. In view of the second identity exhibited preceding the theorem, the corresponding assertion for right alternaive algebas is proved in a similar fashion. Let $ n \in \mathbb{N}$ and  $ f =   \sum_{i=0}^n a_i x^i  \in  \mathbb{A}[x] $ with $ a_n \in \mathbb{A}^{-1}$. Proceed by contradiction and assume that  
$$ f(p) := \sum_{i=0}^n a_i p^i  \in  \mathbb{A}^{-1}, $$ 
for all $ p \in \mathbb{C}_I$.  Thus 
$$ f(x+ I y) = \sum_{i=0}^n a_i (x + I y)^i \in  \mathbb{A}^{-1}, $$
for all $ x , y \in  \mathbb{R}$. Since  $ a_n \in  \mathbb{A}^{-1}$ and $ n \in \mathbb{N}$, we get that 
\begin{eqnarray*}
 \lim_{p \in \mathbb{C}_I,  |p| \rightarrow + \infty}   f(p)^{-1}  & = & \lim_{p \in \mathbb{C}_I,  |p| \rightarrow + \infty} p^{-n} ( a_0 p^{-n}  + \cdots + a_{n-1}p^{-1} + a_n)^{-1}  \\
  & = &  0 \times  a_n^{-1} = 0.
\end{eqnarray*}
Note that by \cite[Theorem 1.21]{Ru} and its real counterpart, the only Hausdorff vector space topology on $ \mathbb{C}_I$  is the Euclidean topology of  $\mathbb{C}_I$. Thus,  $ \lim_{p \in \mathbb{C}_I,  |p| \rightarrow + \infty} p^{-n} =0$.

Let $ \phi \in \mathbb{A}^*:= \mathcal{L}_c(\mathbb{A}, \mathbb{R})$  be an arbitrary continuous real linear functional. We obtain a contradiction by showing that $  \phi( f(p)^{-1} ) = 0$ for all $ p \in \mathbb{C}_I$, implying that $ f(p)^{-1} =0$ by the hypothesis that $ \mathbb{A}^*$ separates the elements of $\mathbb{A}$, a contradiction. Define $ u : \mathbb{R}^2 \longrightarrow \mathbb{R}$ by $ u(x, y) =  \phi( f(x + I y)^{-1} )$. 
Since the  algebra  $ \mathbb{A}$ is alternative, we see that 
\begin{eqnarray*}
 f(x + h+ I y)^{-1}  - f(x +  I y)^{-1}    &= &   \\
  f(x + h+ I y)^{-1} \Big[ \Big( f(x +  I y)  -  f(x + h+ I y) \Big) f(x +  I y)^{-1} \Big],   \\
f\big(x + I( y  + h) \big)^{-1}  - f (x +  I y)^{-1}    & = &   \\
f\big(x + I( y  + h) \big)^{-1} \Big[   \Big(   f\big(x + I( y  + h) \big)  -  f(x + I y) \Big)  f(x + I y)^{-1}\Big],
\end{eqnarray*}

It is easily checked that $ u : \mathbb{R}^2 \longrightarrow \mathbb{R}$ is continuous and that
$$ \frac{\partial^2 u}{\partial x^2} + \frac{\partial^2 u}{\partial y^2} = 0.$$
That is, $u$ is a harmonic function.
We have $ \lim_{||(x, y)||_2 \rightarrow + \infty}  u(x, y) = 0$, because $ \phi$ is continuous, $ \phi(0) = 0$, and that $\lim_{ |p| \rightarrow + \infty}   f(p)^{-1} =0$. From this, we conclude that the harmonic function $ u : \mathbb{R}^2 \longrightarrow \mathbb{R}$  is bounded, and hence it is constant. This implies that $ u = 0$ because $ \lim_{ ||(x, y)||_2 \rightarrow + \infty}  |u(x, y)| = 0$. Therefore, $\phi( f(x + I y)^{-1} ) = 0$ for all $ \phi \in \mathbb{A}^*$ and all $ x, y \in \mathbb{R}$, from which we obtain $f(x + I y)^{-1} = 0$, a contradiction. 

(ii) The proof,  omitted for the sake of brevity, is an imitation of the standard proof of the fundamental theorem of algebra utilizing Liouville's Theorem.
\hfill \qed 

\bigskip

Motivated by the preceding theorems and the main result of  \cite{EN}, here are counterparts of the fundamental theorem of algebra for noncommutative polynomials in the settings of real (resp. complex) left or right  alternative topological algebras.

\bigskip

\begin{thm} \label{2.4} 
  {\rm (i)}  Let   $ \mathbb{A}$ be  a left (resp. right)  alternative  topological real unital algebra whose multiplication is continuous and whose topological dual separates its points, e.g., a locally convex real algebra whose multiplication is jointly continuous. Assume further that   $ \mathbb{A}$ contains a nuclear copy of the complex numbers, namely,  $ \mathbb{C}_I := \{ a + b I : a , b \in \mathbb{R} \}$ with $ I \in {\rm N}(\mathbb{A}) $  and $I^2 = -1$. Then, every nonconstant and noncommutative polynomial with coefficients in $\mathbb{A}$ whose leading noncommutative polynomial part evaluated at some $A_0 \in \mathbb{A}$ is invertible has a singular element  in  $ \mathbb{C}_I A_0$.

{\rm (ii)}   Let  $ \mathbb{A}$ be a left (resp. right)   alternative topological complex unital algebra whose multiplication is continuous and whose topological dual separates its points, e.g.,  a locally convex complex algebra whose multiplication is jointly continuous. Then, every nonconstant and noncommutative polynomial with coefficients in $\mathbb{A}$ whose leading noncommutative polynomial part evaluated at some $A_0 \in \mathbb{A}$ is invertible has a singular element in $ \mathbb{C}A_0$.  
\end{thm}  

\bigskip

\noindent {\bf Proof.} (i)  Just as argued in the proof of the preceding theorem, it suffices to prove the assertion for left alternative algebras. So suppose that the algebra $ \mathbb{A}$ in the statement of the theorem is left alternative. Clearly,   $ \mathbb{A}$ can be viewed as a left alternative  topological  $ \mathbb{C}_I$-algebra.  Note that the complex scalar product $. : \mathbb{C}_I \times  \mathbb{A}  \longrightarrow  \mathbb{A}  $, naturally induced by the multiplication of the real algebra $ \mathbb{A}$,  is continuous because the multiplication of the real algebra  $ \mathbb{A}$ is continuous.  Also note that the complex dual of $ \mathbb{A}$, namely,  $ \mathbb{A}_I^* :=\mathcal{L}_c (  \mathbb{A}, \mathbb{C}_I) $, separates the points of  $ \mathbb{A}$ because so does  the real topological dual  of  $ \mathbb{A}$, namely,  $\mathbb{A}^* := \mathcal{L}_c (  \mathbb{A}, \mathbb{R}) $. To see this, use the hypothesis that $\mathbb{A}^*$ separates the points of $ \mathbb{A}$ and note that   every $ \phi \in  \mathbb{A}^*$ gives rise to a  $ \phi_I \in  \mathbb{A}_I^*$, which is defined by $ \phi_I ( a) = \phi( a) - I \phi( Ia)$. 

 Let $ p = p_n + \cdots + p_1 + p_0$ be a noncommutative polynomial  of degree $n \in \mathbb{N}$, where $p_0 \in \mathbb{A}$ and $ p_i$  ($1 \leq i \leq n$)  is a finite sum of noncommutative monomials of degree $i$ with $p_n (A_0)$   being invertible for some $A_0 \in \mathbb{A}$. Proceeding by way of contradiction, since   $ \mathbb{C}_I \subseteq N(\mathbb{A}) $, we can write 
\begin{eqnarray*}
 \lim_{z \in \mathbb{C}_I,  |z| \rightarrow + \infty}   p(zA_0)^{-1}  & = & \lim_{z \in \mathbb{C}_I,  |z| \rightarrow + \infty} z^{-n} \Big( p_0 z^{-n}  + p_1(A_0) z^{n-1} +\cdots  \\
  & &  \cdots  + p_{n-1}(A_0) z^{-1} + p_n(A_0) \Big)^{-1}  \\
  & = &  0 \times  p_n(A_0)^{-1} = 0.
\end{eqnarray*}
Now, letting $ \phi \in \mathbb{A}_I^*$  be an arbitrary continuous $\mathbb{C}_I$-linear functional, we obtain a contradiction by showing that $  \phi( p(zA_0)^{-1} ) = 0$ for all $ z \in \mathbb{C}_I$, implying that $ p(zA_0)^{-1} =0$ because $  \mathbb{A}_I^*$ separates the points of $ \mathbb{A}$, which is absurd. To this end, define $ f : \mathbb{C}_I \longrightarrow \mathbb{C}_I$ by $ f (z) =  \phi \Big( p \big(zA_0\big)^{-1} \Big)$. Since   $ \mathbb{C}_I \subseteq {\rm N}(\mathbb{A}) $, it is quite straightforward to check that
$$ f'(z) = -\phi \Big( p(zA_0)^{-1} \big[ \big( p_1(A_0) + 2 p_2(A_0) z + \cdots + n p_n(A_0) z^{n-1} \big)  p(zA_0)^{-1} \big] \Big),$$
for all $ z \in \mathbb{C}_I $. Then again, the entire function $f$ is bounded because $ \lim_{z \rightarrow  \infty}  f(z) = 0$, from which, in view of Liouville's Theorem, we get that $f$ is constant, and hence zero. This completes the proof.

(ii) The proof, which is omitted for the sake of brevity, as in (i) is an imitation of the standard proof of the fundamental theorem of algebra utilizing Liouville's Theorem.
\hfill \qed 

\bigskip

\noindent {\bf Remark.}  Motivated by the theorem we suggest the following conjecture. {\it Let   $ \mathbb{A}$ be  a left or right alternative  topological real algebra whose multiplication is continuous and whose topological dual separates its points. Assume further that   $ \mathbb{A}$ contains a copy of the complex numbers, say, 
 $ \mathbb{C}_I := \{ a + b I : a , b \in \mathbb{R} \}$ with $ I \in \mathbb{A} $  and $I^2 = -1$, and $A_0 \in \mathbb{A}$. Then, every nonconstant and noncommutative polynomial with coefficients in $\mathbb{A}$ whose leading noncommutative polynomial part evaluated at every $zA_0 \in \mathbb{A}$ (resp.  $A_0z \in \mathbb{A}$) with $z=  a + b I \in  \mathbb{C}_I $ and $ a^2 + b^2 = 1$ is invertible has a singular element  in  $ \mathbb{C}_I A_0$ (resp. $  A_0 \mathbb{C}_I$).  }

\bigskip

With the preceding theorem at our disposal, the following result is immediate. 

\bigskip

\begin{thm} \label{2.5} 
Let $ n \in \mathbb{N}$
 and  $ \mathbb{A}$ be  an arbitrary  finite-dimensional real or complex unital algebra. If  the algebra $ \mathbb{A}$  happens to be a real algebra, assume further that   $ \mathbb{A}$ contains a nuclear copy of the complex numbers, say,
 $ \mathbb{C}_I := \{ a + b I : a , b \in \mathbb{R} \}$ with $ I \in {\rm N}(\mathbb{A}) $  and $I^2 = -1$. Then, every  element of $ M_n(\mathbb{A})$  has eigenvalues in   $ \mathbb{C}_I$ or $ \mathbb{C}$ depending on whether  $ \mathbb{A}$ is a real or a complex algebra. 

\end{thm}  

\bigskip

\noindent {\bf Proof.} It follows from the hypothesis that $ \mathbb{A}$ is a finite-dimensional algebra over $ \mathbb{C}_I$  or $ \mathbb{C}$ depending on whether  $ \mathbb{A}$ is a real or a complex algebra. So, WLOG, we may assume that $ \mathbb{A}$ is a finite-dimensional complex algebra, and hence it can be normed.  But then again $ \mathbb{A}^n$ can be equipped  with the infinity norm induced by any  algebra norm of  $\mathbb{A}$, and hence $ M_n(\mathbb{A}) \subseteq \mathcal{B} (\mathbb{A}^n)$ can be viewed as an associative complex  normed  algebra together with the operator norm of $ \mathcal{B} (\mathbb{A}^n)$, the complex algebra of all (bounded) linear operators on $\mathbb{A}^n$; note that the elements of $ M_n(\mathbb{A}) $ act on the left of $\mathbb{A}^n$ via the matrix multiplication.  The assertion now follows because, by the preceding theorem,  the noncommutative polynomial $p_A(x) :=  x - A$ with coefficients in $ \mathcal{B} (\mathbb{A}^n)$  has a singular element in   $ \mathbb{C}$ for any given $ A \in  M_n(\mathbb{A})$.  This means that there exists a $ \lambda \in \mathbb{C}$ such that $ \lambda - A$ is singular as a linear operator acting on the left of  $\mathbb{A}^n$, which is finite-dimensional, and hence it has a nonzero kernel. In other words, there exists a nonzero $ X \in \mathbb{A}^n$ such that $ AX = \lambda X = X \lambda$, as desired. 
\hfill \qed 

\bigskip

With a method similar to that used in  the proof of  Theorem \ref{2.2}, we can prove the existence of right eigenvalues for matrices with entries in  finite-dimensional associative real algebras containing a copy of the complex numbers. 

\bigskip

\begin{thm} \label{2.6} 
 Let $ n \in \mathbb{N}$
 and  $ \mathbb{A}$ be  an associative  finite-dimensional real algebra containing a copy of the complex numbers, say, 
 $ \mathbb{C}_I := \{ a + b I : a , b \in \mathbb{R} \}$ with $ I \in \mathbb{A} $  and $I^2 = -1$. Then, every  $ A \in M_n(\mathbb{A})$  has right  eigenvalues in $ \mathbb{C}_I$. In particular, every quaternion matrix has  right  eigenvalues in any copy of the complex numbers within the quaternions. 

\end{thm}  

\bigskip

\noindent {\bf Proof.} By Theorem \ref{1.1} and Proposition \ref{2.1}(i), we can equip the real algebra $ \mathbb{A}$  with an algebra norm having the property that $|| z a w|| = 
 |z| ||a|| |w|$ for all $ z, w \in \mathbb{C}_I $ and $ a \in \mathbb{A}$.  Thus $ \mathbb{A}$ may be viewed as a real normed algebra which is at the same time a complex NLS over  $\mathbb{C}_I $.  We view any $ A   \in M_n(\mathbb{A})$ as a (bounded) linear operator acting on the left of the right  (finite-dimensional) complex normed linear space $ \mathbb{A}^n$ via the matrix multiplication, i.e.,
$A(x) := Ax$. Thus,  $ M_n(\mathbb{A}) \subseteq \mathcal{B} (\mathbb{A}^n)$.  Note that $\mathbb{A}^n$  together with the norm-infinity  $ ||.||_\infty :\mathbb{A}^n \rightarrow \mathbb{R}$, defined by $||x||_\infty = \max ( ||x_1||,  \ldots , ||x_n||)$, where $ x = ( x_1 , \ldots, x_n) \in \mathbb{A}^n $, is a  right complex (finite-dimensional) NLS  over $ \mathbb{C}_I$.   Equip     $ \mathcal{B} (\mathbb{A}^n)$, and hence $ M_n(\mathbb{A})$, with the operator norm induced by the norm-infinity of $\mathbb{A}^n$.  For a given $ p \in   \mathbb{C}_I $, consider the operator $  R_p : \mathbb{A}^n \longrightarrow  \mathbb{A}^n  $, defined by $ R_p(x) = x p$. Clearly, the operators $ R_p$ are bounded linear operators on the (finite-dimensional) complex  NLS $ \mathbb{A}^n$ over $\mathbb{C}_I $. It is plain that $ ||R_p|| = |p|$ for all $ p \in \mathbb{C}_I$, where
$ ||R_p|| := \sup_{||x||_\infty \leq 1}  ||R_p(x)||$. For a $z \in \mathbb{C}_I$ and $ T \in  \mathcal{B} (\mathbb{A}^n)$, recall that $ z T = T z :  \mathbb{A}^n \longrightarrow  \mathbb{A}^n $ is defined by $ (zT)(x) = (Tz)(x) := T(x)z$ ($x \in \mathbb{A}^n $), and hence, we clearly have $ z T = T z \in \mathcal{B} (\mathbb{A}^n)$.

 Let $ A \in M_n(\mathbb{A})$ be given.  Proceed by way of contradiction and suppose that $ A \in M_n(\mathbb{A})$ has no right eigenvalues in $ \mathbb{C}_I$. It follows that the kernel of the complex  bounded linear operator $ A - R_p$ is zero and hence  $ A - R_p$ is invertible for all $ p \in \mathbb{C}_I$  because  $\mathbb{A}^n$ is finite-dimensional.  Define $ f : \mathbb{C}_I \longrightarrow \mathbb{C}_I$ by $ f(x + Iy) =  \phi \Big( \big( A - R_{x + I y}\big)^{-1} \Big)  $, where $ \phi \in \mathcal{B}(\mathbb{A}^n)^* $ is a fixed but arbitrary complex linear functional on $\mathcal{B}(\mathbb{A}^n)$. It is not difficult to check that 
$$f'(x + Iy) = \phi \Big( \big( A - R_{x + I y}\big)^{-2} \Big),$$ 
for all $ z = x + I y \in \mathbb{C}_I$. That is, $ f : \mathbb{C}_I \longrightarrow \mathbb{C}_I$ is an entire function.   
  On the other hand, 
\begin{eqnarray*}
0  \leq |f(x + I y)| & \leq&   || \phi || \ \Big|\Big|  \big( A - R_{x + I y}\big)^{-1}  \Big|\Big|, \\
& = &  || \phi || \  \Big|\Big| R_{x + I y}^{-1} \Big(  AR_{x + I y}^{-1} - I_{\mathbb{A}^n} \Big)^{-1}\Big|\Big|, \\
& \leq& || \phi || \ \Big|\Big| R_{x + I y}^{-1}\Big|\Big|\ \Big|\Big| \Big(  AR_{x + I y}^{-1} - I_{\mathbb{A}^n} \Big)^{-1} \Big|\Big|, \\
& = &  || \phi || \frac{1}{\sqrt{x^2 + y^2}} \Big|\Big| \Big(  AR_{(x + I y)^{-1}} - I_{\mathbb{A}^n} \Big)^{-1} \Big|\Big|,
\end{eqnarray*}
for all $ x , y \in \mathbb{R}$. This clearly yields  $ \lim_{(x, y) \rightarrow  \infty}  |f(x + I y)| = 0$, from which we conclude that the entire function $ f : \mathbb{C}_I \longrightarrow \mathbb{C}_I$  is bounded, and hence it is constant. This implies that $ f = 0$ because  $ \lim_{(x, y) \rightarrow  \infty}  |f(x + I y)| = 0$. Therefore, $\phi\Big( \big( A - R_{x + I y}\big)^{-1} \Big)  = 0$ for all $ \phi \in \mathcal{B}(\mathbb{A}^n)^*$ and all $ x, y \in \mathbb{R}$, from which we obtain $\big( A - R_{x + I y}\big)^{-1}  = 0$, which is absurd. This completes the proof. 
\hfill \qed 

\bigskip


\end{document}